\documentclass[11pt]{amsart}
\usepackage{amsmath,amssymb}
\usepackage[dvips]{graphics}  

\input{epsf}

\theoremstyle{definition}

\theoremstyle{plain}
\newtheorem{theo}{Theorem}[section]
\newtheorem{lem}{Lemma}[section]
\newtheorem{prop}{Proposition}[section]

\newtheorem{cor}{Corollary}[section]

\theoremstyle{remark}
\newtheorem{ps}{Remark}[section]

\def \dR{\mathbb R}
\def \dZ{\mathbb Z}

\def \dC{\mathbb C}
\def \dQ{\mathbb Q}
\def \dS{\mathbb S}
\def \dH{\mathbb H}

\def \sl{SL(2,\dR)}
\def \psl{PSL(2,\dR)}
\def \ch{\text{\rm ch\,}}
\def \sh{\text{\rm sh\,}}
\def \Homeop{\text{Homeo}^+}
\def \Hom{\text{\rm Hom}}
\def \ho{\Hom(\Gamma_g,\psl)}
\def \Tr{\text{\rm Tr}}
\def \tr{\text{\rm tr}}

\def \inject{\hookrightarrow}
\def \vertinj{\cap}
\def \Tau{\mathcal O}

\newcommand{\expon}[1]{\text{\rm exp}\left( #1 \right)}

\title{Non-injective representations of a closed surface group into $\psl$}

\author{Louis Funar}
\address{Institut Fourier BP 74, UMR 5582, Universit\'e  Grenoble I,
38402 Saint-Martin-d'H\`eres Cedex, France}
\email{funar@fourier.ujf-grenoble.fr}
\author{Maxime Wolff}
\address{Institut Fourier BP 74, UMR 5582, Universit\'e  Grenoble I,
38402 Saint-Martin-d'H\`eres Cedex, France}
\email{wolff@fourier.ujf-grenoble.fr}

\begin{document}

\begin{abstract}
  Let $e$ denote the Euler class on the space $\ho$ of representations
of the fundamental group $\Gamma_g$ of the closed surface $\Sigma_g$. Goldman
showed that the connected components of $\ho$ are precisely the inverse
images $e^{-1}(k)$, for $2-2g\leq k\leq 2g-2$, and that the components
of Euler class $2-2g$ and $2g-2$ consist of the injective
representations whose image is a discrete subgroup of $\psl$.
We prove
that non-faithful representations are dense in all the other components.
We show that the image of a discrete representation essentially
determines its Euler class. Moreover, we
show that for every genus and possible corresponding Euler class,
there exist discrete representations.
\end{abstract}

\maketitle

\section{Introduction}

Let $\Sigma_g$ be the closed oriented surface of genus $g\geq 2$.
Let $\Gamma_g$ denote its fundamental group, and $R_g$ the
representation space $\Hom(\Gamma_g,\psl)$. Elements of $R_g$ are
determined by the images of the $2g$ generators of $\Gamma_g$,
subject to the single relation defining $\Gamma_g$. It follows that
$R_g$ has a real algebraic structure (see e.g. \cite{CullerShallen}).
Furthermore, being a subset of $\left(\psl\right)^{2g}$, it is naturally
equipped with a Hausdorff topology.

We can define an invariant $e:R_g\rightarrow\dZ$, called the Euler
class, as an obstruction class or as the index of circle bundles
associated to representations in $R_g$ (see \cite{Milnor,Ghys,Goldman}).

In \cite{Goldman}, which may be considered to be
the starting point of the subject,
Goldman showed that
the connected components of $R_g$ are exactly the fibers $e^{-1}(k)$, for
$2-2g\leq k\leq 2g-2$.
He also proved that $e^{-1}(2g-2)$ and $e^{-1}(2-2g)$ consist of those
injective representations whose image is a discrete subgroup of $\psl$.
Milnor and Wood had previously proved that the inequality $|e(\rho)|\leq 2g-2$
holds for all $\rho\in R_g$ (see \cite{Milnor,Wood}). Goldman
\cite{Gold84} (see also \cite{Hitchin}) showed that
every connected component $e^{-1}(k)$ is a smooth manifold of dimension
$6g-3$, except for the component $e^{-1}(0)$, whose singular points are
the elementary representations.

These connected components have been studied further.
The group $\psl$ acts on $R_g$ by conjugation, and the quotient of
$e^{-1}(2-2g)$ (respectively $e^{-1}(2g-2)$) under this action is the
Teichm\"uller space of $\Sigma_g$ (respectively, the space of marked
hyperbolic metrics with opposite orientation on $\Sigma_g$).

Gallo, Kapovich, Marden and Tan \cite{Galloetal,Tan} showed that
all the representations
(in every connected component of $R_g$)
satisfy a weaker metric condition. Indeed, every representation in $R_g$
is the holonomy representation
of a branched $\dC {\text{\rm P}}^1$-structure on $\Sigma_g$, with
at most one branched point in $\Sigma_g$.
Moreover, Tan showed that for every even $k$ such that $|k|\leq 2g-4$,
there is an explicit non-injective representation in $e^{-1}(k)$, which
is not the holonomy representation of any branched hyperbolic structure on
$\Sigma_g$. However, this explicit representation can be deformed into
representations which are holonomy representations of branched hyperbolic
structures on $\Sigma_g$. It is still unknown whether
the holonomy representations of branched hyperbolic
structures on $\Sigma_g$ form a dense subset of $R_g$.

In this text we take a more elementary point of view. We often consider
representations as products of commutators in $\psl$, and most of our
results involve explicit representations.

Our first result is the following:
\begin{theo}
  For all $g\geq 2$ and all $k$ such that $|k|<2g-2$, non-faithful
  representations form a dense subset of the connected
  component $e^{-1}(k)$.
\end{theo}
Recently Breuillard, Gelander, Souto and Storm
\cite{Suoto}, and independently
DeBlois and Kent IV \cite{DBK} proved that the set of faithful
representations is
also dense, as the intersection of a countable family of
open and dense sets. Therefore theorem 1.1 shows that the
set of non-faithful representations should be thought of as $\dQ$
in $\dR$.

We then show that the Euler class of a discrete representation
is essentially determined by its image, as an abstract group.
  If $\Gamma$ is a non-cocompact Fuchsian group, set $e(\Gamma)=0$.
  Otherwise, if $(g;k_1,\ldots,k_r)$ is the signature of the
  Fuchsian group $\Gamma$
  (with all the $k_i$'s finite), let $d$ be the least common multiple
  of $k_1$, ..., $k_r$ (put $d=1$ if $r=0$) and set
  $$e(\Gamma)=d\left(2g-2+\displaystyle{\sum_{i=1}^{r}}\left(1-\frac{1}{k_i}\right)\right)$$
We have the following result:
\begin{theo}
  Let $\rho$ be a representation (in any $R_{g'}$) whose image is contained in
  $\Gamma$. Then $e(\rho)$ is a multiple of $e(\Gamma)$.
  Moreover, every multiple of $e(\Gamma)$ is the Euler class of some
  representation whose image is exactly $\Gamma$.
\end{theo}
\noindent In particular, there are no representations of non-zero Euler
class taking values in $PSL(2,\dZ)$.

\noindent We deduce the following proposition:
\begin{prop}
  Let $k$ be a fixed non-zero integer. There are only finitely many Fuchsian
  groups $\Gamma$ such that there exists a representation (of $\Gamma_g$, for
  some $g\geq 2$)
  with image contained in $\Gamma$ and with Euler class $k$.
\end{prop}
\noindent Moreover, we have:
\begin{prop}
  Let $g\geq 2$ and $k\in\dZ$ be such that $|k|\leq 2g-3$, $k\neq 0$. The
  representations whose image is discrete form a nowhere dense closed subset
  of $e^{-1}(k)$ in $R_g$.
\end{prop}
\noindent In other words, there are few discrete representations of non-zero
Euler class.
On the other hand,
we have the following:
\begin{theo}
  For every $g\geq 2$ and every $k$ such that $|k|\leq 2g-2$,
  there exists a discrete representation of Euler class $k$ in $R_g$.
\end{theo}
\noindent Those representations are given explicitly in terms of
signatures of Fuchsian groups.
In particular, using Magnus' results from \cite{Magnus}, these representations
can be expressed using only matrices in
$PSL(2,\dZ\left[\frac{1}{2}\right])$.

Finally, we deduce a characterization of representations of odd Euler
class, which enables us to describe all the subgroups of $\psl$
that are the image of some discrete representation of odd Euler
class (in some $R_g$). If $\Gamma$ is a cocompact Fuchsian
group with signature
$(g';k_1,k_2,\ldots,k_r)$, let $m(\Gamma)$ be the maximal power
of $2$ dividing one of the $k_i$'s. If $m(\Gamma)=0$, set
$n(\Gamma)=0$. Otherwise, let $n(\Gamma)$ be the number
of $k_i$'s which are divisible by $2^{m(\Gamma)}$.
We have  then the following characterization:
\begin{prop}
  A Fuchsian group $\Gamma \subset \psl$ is the image of a
  representation (of $\Gamma_g$, for some $g\geq 2$)
  if and only if $\Gamma$ is
  a cocompact Fuchsian group such that $n(\Gamma)$ is odd.
\end{prop}

{\bf Acknowledgements.}
The authors are grateful to R. Bacher, M. Eisermann, V. Fock, E. Ghys,
C. MacLean, G. McShane and V. Sergiescu
for inspiring  discussions, comments and corrections.

\section{The Lie group $\psl$ and Milnor's algorithm}

\subsection{The Lie group $\psl$}{\ }

We recall first some basic results on $\psl$ and we refer to \cite{Katok}
for a full treatment of this subject.

The Lie group $\sl$ is $\{M\in GL(2,\dR) | \det M=1\}$ and 
$\psl$ is the quotient $\sl / \{\pm 1\}$ by its center. 
Topologically, $\sl$ and $\psl$
are two three-dimensional solid tori, and the projection map
$\sl\rightarrow \psl$ is a $2$-sheeted cover. We will denote
elements of $\psl$ by matrices, as if they were in $\sl$.

If $M=\left(\begin{array}{cc}a & b \\ c & d \end{array}\right)\in \psl$,
the homography $z\mapsto\dfrac{az+b}{cz+d}$ is an isometry of the
upper half-plane model of the hyperbolic plane $\dH^2$.
This isometry of $\dH^2$ acts on its boundary
$\dS^1=\dR\cup\{\infty\}$, still by homography. This defines an injection
$\psl\inject\Homeop(\dS^1)$, which we will think of as an inclusion.

\smallskip

On $\psl$ only the absolute value of the trace is well-defined.
If $M\in\sl$, we denote
$\Tr(M)=|\tr(M)|$.
Elements of $\psl$ then fall into three types, $Ell$, $Par$ and $Hyp$,
depending on their traces:

\begin{itemize}
  \item If $\Tr(M)\in[0,2)$ then $M$ fixes a point in $\dH^2$ and
    acts as a rotation around it.
    The matrix $M$ is conjugate {\it in $\psl$} to a matrix of the
    form $\left(\begin{array}{cc}\cos\theta & -\sin\theta \\
    \sin\theta & \cos\theta \end{array}\right)$. It follows that two such
    elements are conjugate if and only if they have the same
    trace (in absolute value). Such matrices are called {\em elliptic}.
  \item If $\Tr(M)=2$ and $M\neq I_2$, then $M$ fixes no point in
    $\dH^2$ but fixes a unique point in its boundary. Its orbits in
    $\dH^2$ are the horospheres defined by that fixed point. The matrix
    is conjugate, in $\psl$, to a matrix
    $\left(\begin{array}{cc}1 & t \\ 0 & 1 \end{array}\right)$. Such
    elements are called {\em parabolic}. Elements
    of $Par$ fall into two conjugacy classes (over $\psl$), depending
    on the sign of $t$.
  \item If $\Tr(M)\in(2,+\infty)$ then $M$ is conjugate in $\psl$ to
    $\left(\begin{array}{cc}\ch u & \sh u \\ \sh u & \ch u \end{array}\right)$
    for some $u\in\dR$,
    and two such elements are conjugate if and only if they have
    the same trace. The action of $M$ on $\dH^2$ fixes
    two points in the
    boundary of $\dH^2$. Those elements are called {\em hyperbolic}.
\end{itemize}

Each of these three sets contains (together with $I_2$) whole one-parameter
groups, which are determined by the fixed points described above.

\smallskip

\subsection{Milnor's algorithm}{\ }

Now we describe Milnor's algorithm for calculating the Euler
class of a representation.
We refer to Ghys' article \cite{Ghys} for a complete and detailed approach
to this important way of computing (or defining) Euler classes.

Topologically, the Lie group $\psl$ is a solid torus, and $\pi_1(\psl)\cong\dZ$.
Moreover, for all $A\in \psl$, taking a lift
$\tilde{A}\in\widetilde{\psl}$ is the same as lifting the homeomorphism
$A\in\Homeop(\dS^1)$ to a homeomorphism of the universal cover of the circle.
%$\widetilde{\dS^1}=\dR$.
In other words, we have the following short exact sequence:
$$ \begin{array}{rcccccccl}
0 & \rightarrow & \dZ & \rightarrow &
\widetilde{\psl} & \rightarrow & \psl & \rightarrow & 1 \\
 & & & & \vertinj & & \vertinj & & \\
 & & & & \Homeop(\dR) & & \Homeop(\dS^1) & &
\end{array}$$

In the diagram above, the sign of the generator in $\dZ$ is determined
by the choice of an orientation of $\dS^1$.
As in \cite{Goldman}, we will denote by $z\in\widetilde{\psl}$ the
image of this generator $1\in\dZ$.

Note that, for all $A\in \psl$, the choice of a lift
$\tilde{A}\in \widetilde{\psl}$ is well-defined up to a certain number of
elementary translations $z$, which are central in $\widetilde{\psl}$.

The Euler class $e(\rho)$ of a representation $\rho$ is 
computed in the following way. 
Consider the standard one relator  presentation of $\Gamma_g$ given 
by 
\[ \Gamma_g =\langle \; a_1,a_2,\ldots,a_g,b_1,b_2\ldots,b_g \;|\;\; [a_1,b_1][a_2,b_2]\cdots [a_g,b_g]\;
\rangle\]
For each generator $x$ in the above presentation of $\Gamma_g$,
choose an arbitrary lift $\widetilde{\rho(x)}\in\widetilde{\psl}$. Then
$e(\rho)$ is determined by means of the formula:
$$[\widetilde{\rho(a_1)},\widetilde{\rho(b_1)}]\cdots
[\widetilde{\rho(a_g)},\widetilde{\rho(b_g)}]=z^{e(\rho)}$$
\noindent Obviously, $e(\rho)$ does not depend on the choice of the 
lifts. Indeed, two lifts
differ by powers of $z$, which disappear in the commutators, because $z$
is central in $\widetilde{\psl}$.

\section{Non-injective representations}
In this section we give first an elementary proof of the
following proposition, which is a particular case of theorem 1.1:

\begin{prop}
  For all $k$ such that $|k|\leq 2g-3$, $e^{-1}(k)$ contains an explicit open
set $E$, in which non-injective representations are dense.
\end{prop}

We set $E=\{\rho\in R_g|\rho(a_1)\in Ell,\rho(b_1)\in Ell\}\cap e^{-1}(k)$.

In order to prove the proposition 3.1,
we first need some preliminaries. We then 
state two lemmas implying the proposition, and finally prove the
two lemmas.

\subsection{Preliminaries}

\begin{lem}
  Every hyperbolic element in $\psl$ is the commutator of two
elliptic elements.
\end{lem}

\begin{proof}
Let $A_\theta=\left(\begin{array}{cc}\cos\theta & -\sin\theta \\ \sin\theta &
\cos\theta \end{array}\right)$ and $U_t=\left(\begin{array}{cc}1 & t \\
0 & 1 \end{array}\right)$. Then $\Tr([A_\theta,U_t A_\theta U_{-t}])=
2+(4t^2+t^4)\sin^4\theta$. This trace takes every value from the 
interval $[2,+\infty)$,
and the set $[Ell,Ell]$ is invariant under conjugation.
Thus every
hyperbolic element is the commutator of two elliptic elements.
\end{proof}

\begin{ps}
  If $M\in Hyp\cup Par$, or equivalently, if $M\in\Homeop(\dS^1)$ has
a fixed point in $\dS^1$,
then there is a canonical lift $\widetilde{M}^{\text{\rm can}}$ of $M$,
in $\Homeop(\dR)$, which has fixed points in $\dR$.
\end{ps}

\begin{ps}
  If $A$, $B\in \psl$, $\widetilde{A}$ and $\widetilde{B}$
are defined only up to a power of the elementary translation
$z\in\Homeop(\dR)$, but the commutator 
$[\widetilde{A},\widetilde{B}]\in \widetilde{\psl}$
is well-defined and independent on the choice of lifts, because 
$z$ is central in $\widetilde{\psl}$.
\end{ps}

\begin{prop}
  For all $\varepsilon\in\{-1,0,1\}$ and for all $M\in Hyp$, there
exist $A,B\in \psl$ such that
$[\widetilde{A},\widetilde{B}]=\widetilde{M}^{\text{\rm can}}\cdot
z^\varepsilon$.
  Moreover, if $[A,B]=M$ then the only possible values of
$[\widetilde{A},\widetilde{B}]^{-1}\widetilde{M}^{\text{\rm can}}$
are $z^{-1}$, $1$ and $z$.
\end{prop}
\noindent Proofs of this important result can be found in
\cite{Milnor,Goldman}.

\noindent Another well-known result is the following:
\begin{prop}
  If $M\in Hyp$, $A,B\in Ell$ and $[A,B]=M$ then
$[\widetilde{A},\widetilde{B}]=\widetilde{M}^{\text{\rm can}}$.
\end{prop}

\begin{proof}
  This is a corollary of Goldman's main result in \cite{Goldman}.
  Indeed, suppose
$[\widetilde{A},\widetilde{B}]\neq\widetilde{M}^{\text{\rm can}}$.
Then
$[\widetilde{A},\widetilde{B}]=
\widetilde{M}^{\text{\rm can}}\cdot z^\delta$,
with $\delta\in\{-1,1\}$.
  Now, $M^{-1}\in Hyp$, so, by proposition 3.2, there exist
$C,D\in \psl$ such that
$[\widetilde{C},\widetilde{D}]=\widetilde{M^{-1}}^{\text{\rm can}}\cdot
z^\delta$.
  Hence
$[\widetilde{A},\widetilde{B}][\widetilde{C},\widetilde{D}]=z^{2\delta}$,
so the formulas 
$$\rho(a_1)=A, \rho(b_1)=B, \rho(a_2)=C, \rho(b_2)=D$$
define 
a representation $\rho$ of the fundamental group
$\Gamma_2=\left\langle a_1,b_1,a_2,b_2 \left|[a_1,b_1][a_2,b_2]
\right.\right\rangle$
of the genus two surface $\Sigma_2$, of
Euler class $2\delta$. Therefore, by  corollary C of 
\cite{Goldman} 
this representation is faithful and discrete, which contradicts
the assumption $A\in Ell$.
\end{proof}

\subsection{Proof of proposition 3.1}

\begin{lem}
  For all $\rho\in E$, there exists a path $\rho_t\in R_g$
such that $\rho_0=\rho$
and $\dfrac{d}{dt}\Tr(\rho_t(a_1))\neq 0$.
\end{lem}

\begin{ps}
  We are not concerned with the  derivability of the
  function $t\mapsto\Tr(\rho_t(a_1))$, but  rather in  the fact that 
  it is not locally constant.
\end{ps}

\begin{lem}
  The set $E$ is not empty.
\end{lem}

An element $A\in Ell$ is conjugate to a unique matrix
$\left(\begin{array}{cc}\cos\theta & -\sin\theta \\
\sin\theta & \cos\theta \end{array}\right)$.
If $\dfrac{d}{dt}\Tr(\rho_t(a_1))$ is non-zero, then the angle
$\theta$ corresponding to $\rho_t(a_1)$ is rational for a dense
subset of the possible values of $t$, hence the representation
$\rho_t$ satisfies a relation of the kind $\rho_t(a_1)^n=1$, so
that $\rho_t$ is not injective.
Thus these two lemmas imply proposition 3.1.

The two remaining sections consist of proofs of these two
lemmas.

\subsection{Proof of lemma 3.2}{\ }

Let $A=\rho(a_1)$ and $B=\rho(b_1)$. Let $B(t)$ be the one-parameter
subgroup of $PSL_2(\dR)$ containing $B$. That is, $B(0)=1$, $B(1)=B$
and for all $t\in\dR$, $B(t)$ commutes with $B$.
Set $A_t=AB(t)$. Then $[A_t,B]=AB(t)BB(-t)A^{-1}B^{-1}=[A,B]$. 
Thus the formulas  
\[ \rho_t(a_1)=A_t,\; \rho_t(a_i)=\rho(a_i), \;{\rm for } \;i\geq 2,\; 
\rho_t(b_j)=\rho(b_j), \;{\rm for } \;j\geq 1 \]
define a representation $\rho_t\in\ho$. Moreover, $\rho_0=\rho$.

We claim that $\dfrac{d}{dt}\Tr(\rho_t(a_1))\neq 0$.
Up to conjugation -- which does not change the value of the
traces -- we can write $A$ in an adapted basis under the form
$A=\left(\begin{array}{cc}\cos\theta & -\sin\theta \\
\sin\theta & \cos\theta \end{array}\right)$. In the same basis $B(t)$
reads 
$\left(\begin{array}{cc}c_1(t) & c_2(t) \\ c_3(t) & c_4(t)
\end{array}\right)$,
where $c_j(t)$, $1\leq j\leq 4$, are  suitable real functions. 
Set $d_i=c_i'(0)$.
We compute 
\[  {\rm tr}(AB(t))=(c_1(t)+c_4(t))\cos\theta+(c_2(t)-c_3(t))\sin\theta\]
and hence 
\[ \dfrac{d}{dt}{\rm tr}(AB(t))=(d_1+d_4)\cos\theta+(d_2-d_3)\sin\theta\]
It suffices to show that this derivative is non-zero.  
Notice that $d_1+d_4=0$, and
$A\in Ell$ (which does not contain $1$) so $\sin\theta\neq 0$.
Hence, the claim  would follow if we  proved that  $d_2\neq d_3$. 
Assume  the contrary, namely that $d_2=d_3$. 
We develop $B(t)$ up to second order terms:
$$B(t)=\left(\begin{array}{cc}1+d_1 t+\alpha t^2+\Tau(t^3) & d_2 t+\Tau(t^2) \\
d_2 t+\Tau(t^2) & 1-d_1 t+\beta t^2+\Tau(t^3)\end{array}\right) $$
Now, $\det B(t)=1-(\alpha+\beta-d_1^2-d_2^2)t^2+\Tau(t^3)=1$, which implies 
that $\alpha+\beta>0$. Therefore, for small enough $t>0$ (and thus for all
$t$), $B(t)$ is hyperbolic, which is a contradiction since
we assumed that $B\in Ell$.

\subsection{Proof of lemma 3.3}{\ }

If $|k|\leq 2g-4$, take a representation $\rho'\in R_{g-1}$  
of Euler class $k$. Assume that the standard generators of $\Gamma_{g-1}$ 
are denoted by the letters  $a_1',b_1',\ldots a_{g-1}',b_{g-1}'$ in order 
to avoid confusion with the generators $a_1,a_2,\ldots,a_g,b_g$ of 
$\Gamma_g$. 
%\[\Gamma_{g-1}=\left\langle a_1',b_1',\ldots a_{g-1}',b_{g-1}'\left|
%[a_1',b_1']\cdots[a_{g-1}',b_{g-1}']\right.\right\rangle\]
Choose an arbitrary $A\in Ell$.
Set $\rho$ for the homomorphism defined by the formulas 
\[ \rho(a_1)=\rho(b_1)=A, \rho(a_i)=\rho'(a_{i-1}'), \; \rho(b_i)=
\rho'(b_{i-1}'), \; {\rm for}\; i\geq 2 \] 
Then $\rho\in E$.

If $|k|=2g-3$, it suffices to consider the case $k=2g-3$.
Then take a representation
$\rho'\in R_g$ of Euler class $2g-2$.
We can suppose (see e.g. \cite{Goldman})
that $\rho'([a_1,b_1])$ is hyperbolic, so by lemma 3.1
there exist $A,B\in Ell$ such that $[A,B]=\rho'([a_1,b_1])$.
The formulas 
$$\rho(a_1)=A\text{, }\rho(b_1)=B\text{ and }
\rho(a_i)=\rho'(a_i)\text{, }\rho(b_i)=\rho'(b_i),\text{ for } \;2\leq i\leq g$$
define a
representation $\rho\in\ho$.
In order to prove that $\rho\in E$,
we need to check that $e(\rho)=k$.
But by maximality of the Euler class of $\rho'$ and proposition 3.2,
we know that $[\widetilde{\rho'(a_1)},\widetilde{\rho'(b_1)}]=
\widetilde{\rho'([a_1,b_1])}^{\text{\rm can}}\cdot z$, and by proposition
3.3 we know that $[\widetilde{A},\widetilde{B}]=
\widetilde{\rho'([a_1,b_1])}^{\text{\rm can}}$. Thus 
$e(\rho)=e(\rho')-1=k$, as claimed.

\bigskip

Our proof of theorem 1.1 uses theorem 1.3. We postpone the proof of 
both theorems for  the next section.

\section{Representations with discrete image}

\subsection{Presenting lifts of Fuchsian groups into $\widetilde\psl$}{\ }

Our main technical result, which will be used throughout this sequel,
is the following:
\begin{lem}
  Let $\Gamma$ be the Fuchsian group with signature
  $(g;k_1,\ldots,k_l,\underbrace{\infty,\ldots,\infty}_{r-l})$.
  Then the lift $\widetilde{\Gamma}$ of $\Gamma$ in $\widetilde\psl$
  has the following presentation:
  $$\widetilde{\Gamma}=\left\langle q_1,\ldots,q_r,a_1,\ldots,b_g,z\left|
  \begin{array}{l}zq_1z^{-1}q_1^{-1},
  \ldots,
  zb_gz^{-1}b_g^{-1}, \\ q_1^{k_1}z,\ldots,q_l^{k_l}z,
  q_1\cdots q_r[a_1,b_1]\cdots[a_g,b_g]z^{2g-2+r}
  \end{array}\right.\right\rangle$$
\end{lem}

\begin{proof}
  Take a Poincar\'e  fundamental
  domain of $\Gamma$ in  the hyperbolic plane $\dH^2$. 
     According to Katok's notations
  (\cite{Katok}), we have the following picture in the example case
  where $g=2$, $r=4$, $k_i>2$ for $i=1,2,3$ and $k_4=2$:

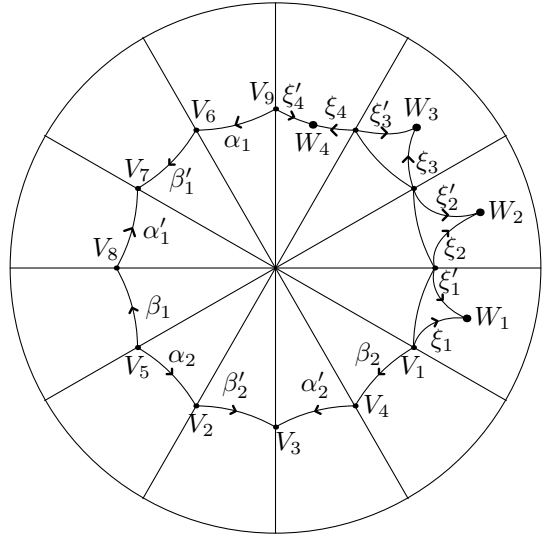
\begin{center}
\includegraphics{funddom1.eps}
\end{center}  
%$$\epsfbox{funddom.1}$$
 
 Here $q_i$ is the only element of $\psl$ sending $\xi_i'$ to
  $\xi_i$, and similarly $a_i(\alpha_i')=\alpha_i$ and
  $b_i(\beta_i')=\beta_i$.
  
  Now consider the identification
  of $\psl$ with the unit tangent bundle
  $S\dH^2$ of the hyperbolic
  plane $\dH^2$ (see \cite{Katok}, theorem 2.1.1).
  Then the universal lift $\widetilde{\psl}$ of $\psl$
  can be viewed  as the unit tangent bundle $S\dH^2$ together
  with an index number. And lifting elements of $\psl$ to $\sl$
  is just taking the index number modulo 2.
  Given our prefered presentation of $\Gamma$, in order to determine
  a presentation of $\widetilde{\Gamma}$ we just have to know how
  the relations lift. One may lift $q_i$ to the rotation
  of positive angle $\frac{2\pi}{k_i}$, so that $q_i^{k_i}$ makes
  one (positive) turn around $W_i$. Hence the relation $q_i^{k_i}$
  lifts to $q_i^{k_i}z^{-1}$ in $\widetilde{\psl}$. We now have to
  determine the lift of the long relation
  $q_1\cdots q_r[a_1,b_1]\cdots[a_g,b_g]$
  in $\widetilde{\Gamma}$.
  
  For this, we consider two vectors $s_i$ and $s_i'$ based at
  each point $V_i$ or $W_i$.

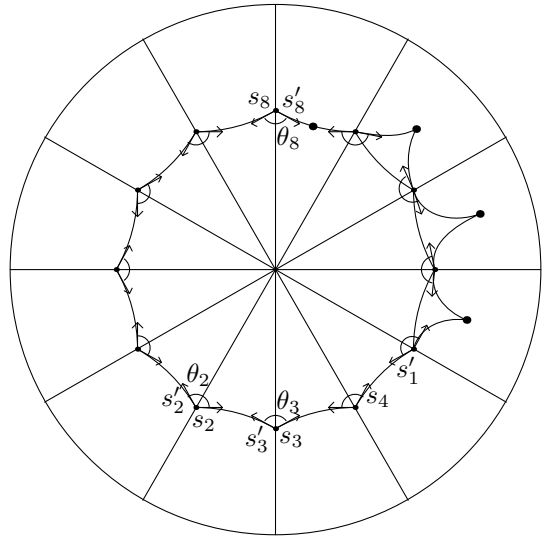
\begin{center}  
\includegraphics{funddom2.eps}
\end{center}
%$$\epsfbox{funddom.2}$$
 
 The element $a_i$ is defined (see \cite{Katok}) as the only element
  of $\psl$ sending $\alpha_i'$ to $\alpha_i$; we choose its lift in
  $\widetilde{\psl}$ to be the transformation of $S\dH^2$
  sending $s_{4g-4i+3}$ to $s_{4g-4i+2}'$. The same element also sends
  $s_{4g-4i+4}'$ to $s_{4g-4i+5}$. Similarly, for the lift of $b_i$ to
  $\widetilde{\psl}$, take the transformation of $S\dH^2$ sending
  $s_{4g-4i+3}'$ to $s_{4g-4i+4}$ and $s_{4g-4i+2}$ to $s_{4g-4i+1}'$.
  Then $b_g^{-1}(s_1')=s_2$, $a_g^{-1}(s_2')=s_3$, ..., $q_1(s_{4g+r}')=s_1$.
  Now define
  $\Theta_i: S\dH^2 \rightarrow S\dH^2$
  by $\Theta_i(p,v)=(p,v+\theta_i)$.
  Then $\Theta_i(s_i)=s_i'$, and
  $\Theta_i$ is central in $\widetilde{\psl}$. Moreover,
  $\Theta_{4g+r}\circ\cdots\circ\Theta_1=z$, or equivalently,
  $\displaystyle{\sum_{i=1}^{4g+r}}\theta_i=2\pi$ (by the Gauss-Bonnet
  formula (see also \cite{Katok}, theorem 4.3.2)). Now, if
  the element $q_1\cdots q_r[a_1,b_1]\cdots[a_g,b_g]$
  caused, say, $n$ turns in the index
  (i.e. $q_1\cdots q_r[a_1,b_1]\cdots[a_g,b_g]=z^n$),
  then we make $n+1$ turns after following the arrows $s_1'$, $s_2$,
  $s_2'$,..., $s_{4g+r}'$, $s_1$ and back to $s_1'$. Now we need
  just compute this index in order to find $n$.
  
  The simplest way to do this is to extend our arrows to a
  vector field whose only singular points
  in the surface $\dH^2/\Gamma$ are the central point $O$ in the
  center of our fundamental domain, and the point $W$ coming from the
  $V_i$'s and $W_i$'s. The index of this vector field at the point $O$ is
  $-(n+1)$ (the arrows point towards $O$, hence the minus sign).
  And since $\displaystyle{\sum_{i=1}^{4g+r}}\theta_i=2\pi$,
  the index at $W$ is $1$. Hence, by the classical
  Poincar\'e-Hopf index theorem,
  $$-(n+1)+1=\chi(\Sigma_g^r)=2g-2+r\text{\rm, i.e. }n=2g-2+r$$
  This gives the
  following presentation for $\widetilde{\Gamma}$:
  $$\widetilde{\Gamma}=\left\langle q_1,\ldots,q_r,a_1,\ldots,b_g,z\left|
  \begin{array}{l}zq_1z^{-1}q_1^{-1},
  \ldots,
  zb_gz^{-1}b_g^{-1}, \\ q_1^{k_1}z^{-1},\ldots,q_l^{k_l}z^{-1},
  q_1\cdots q_r[a_1,b_1]\cdots[a_g,b_g]z^{2-2g-r}
  \end{array}\right.\right\rangle$$
  which is obviously equivalent to our statement (replace $z^{-1}$
  with $z$).
\end{proof}

\begin{ps}
  We could also count this index geometrically.
  First suppose for simplicity that
  $g=0$.  Then the tail of our vector pointed to the center of the fundamental
  domain exactly once for every rotation, thus $r$ times. And we also
  turned around this center once negatively, so the index is
  $r-1$. Thus $n=r-2$, or equivalently, the element
  $q_1\cdots q_r[a_1,b_1]\cdots[a_g,b_g]$ lifts to
  $q_1\cdots q_r[a_1,b_1]\cdots[a_g,b_g]h^{-r+2}$ in
  $\widetilde{\psl}$.
  In the general case ($g\geq 1$) one notes that
  each of the commutators
  $[a_i,b_i]$ gives two positive turns to our arrow.
\end{ps}

\subsection{The Euler class is determined by the image, for discrete 
representations}{\ }

We will present here the proof of  theorem 1.2.
  Consider a Fuchsian group $\Gamma$ with signature
  $(g;k_1,\ldots,k_l,\underbrace{\infty,\ldots,\infty}_{r-l})$.
  If $r\neq l$,
  set $e(\Gamma)=0$. Otherwise,
  let $d$ be the least common multiple
  of $k_1$, ..., $k_r$ (put $d=1$ if $r=0$) and set
  $$e(\Gamma)=d\left(2g-2+\displaystyle{\sum_{i=1}^{r}}\left(1-\frac{1}{k_i}\right)\right)$$
\begin{prop}
  Let $\rho$ be a representation of some $\Gamma_{g'}$ 
  whose image is contained in
  $\Gamma$. Then $e(\rho)$ is a multiple of $e(\Gamma)$.
  Moreover, every multiple of $e(\Gamma)$ is the Euler class of some
  representation whose image is exactly $\Gamma$.
\end{prop}

\begin{proof}
  It follows from lemma 4.1 that $H^1(\widetilde{\Gamma})$ has the
  following abelian presentation:
  $$H^1(\widetilde{\Gamma})=\left\langle q_1,\ldots,q_r,a_1,\ldots,b_g,z\left|
  q_1^{k_1}z,\ldots,
  q_l^{k_l}z,q_1\cdots q_r z^{2g-2+r}\right.\right\rangle^{\rm ab}$$
where all generators are supposed to commute with each other, 
and this is recorded by using the superscript $ab$ on the right side. 

  We just need to show that the subgroup $\langle z\rangle$ generated by $z$
  in $H^1(\widetilde{\Gamma})$ is exactly $\dZ/e(\Gamma)\dZ$, where
  $e(\Gamma)=d\left(2g-2+\displaystyle{\sum_{i=1}^{r}}\left(1-\frac{1}{k_i}\right)\right)$.
  Now $z^N=1$ in $H^1(\widetilde{\Gamma})$ if and only if $z^N$ is a
  product of relations $q_1^{k_1}z$,...,$q_l^{k_l}z$ and
  $q_1\cdots q_r z^{2g-2+r}$ in the {\em abelian} group
  $H^1(\widetilde{\Gamma})$.
  
  Let us write these relations as column vectors in terms of the
  generators $q_1,...,q_r$ and $z$ (in the case $l=r$):
  $$\left(\begin{array}{ccccc}k_1 & 0 & \ldots & 0 & 1 \\
  0 & k_2 & \ddots & \vdots & 1 \\ \vdots & \ddots & & 0 & \\
  0 & \ldots & 0 & k_r & 1 \\ 1 & 1 & \ldots & 1 & 2g-2+r \end{array}\right)$$
  Now, suppose  that $z^N=1$, i.e. there exists a linear relation 
  $$\left(\begin{array}{c} 0 \\ 0 \\ \vdots \\ 0 \\ N \end{array}\right)=
  \alpha_1\left(\begin{array}{c}
  k_1 \\ 0 \\ \vdots \\ 0 \\ 1 \end{array}\right)+
  \alpha_2\left(\begin{array}{c}
  0 \\ k_2 \\ 0 \\ \vdots \\ 1 \end{array}\right)+ \cdots +
  \alpha_r\left(\begin{array}{c}
  0 \\ \vdots \\ 0 \\ k_r \\ 1 \end{array}\right)+
  \beta\left(\begin{array}{c} 1 \\ 1 \\ \vdots \\ 1 \\ 2g-2+r
  \end{array}\right)$$
  with $\alpha_1$,..., $\alpha_r$, $\beta$ integers. Then,
  $k_1\alpha_1+\beta=0$,..., $k_r\alpha_r+\beta=0$.
  Hence, $\beta$ is a multiple of each of the $k_i$'s, i.e.
  $\beta=nd$ for some $n\in\dZ$. Therefore $\alpha_i=\frac{-nd}{k_i}$.
  Now, $$N=\alpha_1+\cdots+\alpha_r+\beta(2g-2+r)=
  nd\left(2g-2+r-
  \displaystyle{\sum_{i=1}^{r}}\frac{1}{k_i}\right)=
  nd\left(2g-2+\displaystyle{\sum_{i=1}^{r}}
  \left(1-\frac{1}{k_i}\right)\right)$$
  Conversely, we can write $z^{ne(\Gamma)}$ explicitly as a product
  of commutators, as it will be done in the next remark. This defines a
  representation of Euler class $ne(\Gamma)$ having its image contained in
  $\Gamma$. Now, to find a representation whose image is exactly $\Gamma$,
  just multiply $z^{ne(\Gamma)}$ by
  $[q_1,q_1]\cdots[q_r,q_r][a_1,a_1]\cdots[b_g,b_g]$
  (as in remark 4.7). This yields $z^{ne(\Gamma)}$ as a
  product of commutators, thus defining a representation of the group
  of some surface $\Sigma_{g'}$, of Euler class $ne(\Gamma)$.
  
  And in the case $l<r$, the $r^{\rm th}$ line simply gives $\beta=0$, so
  that $N=0$, which
  completes the proof.
\end{proof}

\begin{ps}
  Given a Fuchsian group $\Gamma$, let $G(\Gamma,n)$ be
  the minimal genus such that there exists a representation
  $\rho\in R_{G(\Gamma)}$ of image (contained in) $\Gamma$ and
  Euler class $n e(\Gamma)$.
  The Milnor-Wood inequality together with explicit formulas
  give the rough inequalities
  $$n\frac{e(\Gamma)}{2}+1\leq G(\Gamma,n)\leq ndg+r^{nd}$$
  Indeed, write $z^{ne(\Gamma)}$ as a product of $G(\Gamma,n)$
  commutators in $\widetilde{\Gamma}$: then
  $ne(\Gamma)\leq 2G(\Gamma,n)-2$, i.e.
  $G(\Gamma,n)\geq n\frac{e(\Gamma)}{2}+1$. Concerning  the upper
  bound, one can check by induction on $d$ that the product
  $(q_1\cdots q_r)q_1^{-d}\cdots q_r^{-d}$
  can be written as a product of $r^d$ commutators. Now in
  $\widetilde{\Gamma}$,
  $\left(q_1\cdots q_r c_g z^{2g-2+r}\right)^{nd}=1$
  so $\left(q_1\cdots q_r c_g\right)^{nd}q_1^{-nd}
  \cdots q_r^{-nd}=z^{-ne(\Gamma)}$ and
  $\left(q_1\cdots q_r c_g\right)^{nd}q_1^{-nd}
  \cdots q_r^{-nd}$ can be written as the product
  of $ndg+r^{nd}$ commutators.
  
  It is quite easy to refine these inequalities; however, we do
  not have sharp estimates for $G(\Gamma,n)$.
\end{ps}

\begin{ps}
  The term $2g-2+\displaystyle{\sum_{i=1}^r}\left(1-\frac{1}{k_i}\right)$
  is the volume of the Poincar\'e  fundamental domain of the Fuchsian group
  $\Gamma$ (see e.g. \cite{Katok}). In particular, it is positive; hence
  $e(\Gamma)$ is a positive integer if $\Gamma$ is a cocompact Fuchsian
  group.
\end{ps}

\begin{ps}
  Of course, $l\neq r$ if and only if $\Gamma$ is non-cocompact. In this
  case, there is another reason for every representation taking
  values in $\Gamma$ to have zero Euler class. If $\Gamma$ is
  non-cocompact, then $H^2(\Gamma)=0$ (since the surface
  $\dH^2 / \Gamma$ is not compact). And we have the following:
\end{ps}

\subsection{Restrictions on Fuchsian image groups for fixed Euler class}{\ }

\begin{prop}
  Let $\Gamma$ be any group injectively
  mapped in $\psl$ such that $H^2(\Gamma)=0$,
  and let $\rho$ be a representation with image contained in $\Gamma$.
  Then $e(\rho)=0$.
\end{prop}

\begin{proof}
  The short exact sequence
  $0\rightarrow\dZ\rightarrow\widetilde{\Gamma}\rightarrow\Gamma\rightarrow 0$
  is a central extension so the
  corresponding spectral extension is exact at its end (see e.g. \cite{Brown}):
  $$\cdots\rightarrow H^2(\Gamma)\rightarrow H^1(\dZ)\rightarrow
  H^1(\widetilde{\Gamma})\rightarrow \cdots$$
  Now, $H^1(\dZ)=\dZ$
  and the image of any non zero $n\in H^1(\dZ)$, is not zero since
  $H^2(\Gamma)=0$. In other words, in $\widetilde{\Gamma}$
  no power of $z$ is a product of commutators.
\end{proof}

\begin{ps}
  However, there exist groups $\Gamma$, embedded in $\psl$, with
  $H^2(\Gamma)\neq 0$ such that every representation
  taking values in $\Gamma$ has Euler class $0$. For instance, every
  realization of $\dZ^2$ in $\psl$ is in a one-parameter subgroup of
  $\psl$, and hence every representation taking values in
  $\dZ^2\subset\psl$ is elementary
  and thus has Euler class $0$.
\end{ps}

On the other hand, many groups are the image of some representation
of zero Euler class:
\begin{ps}
  For every
  finitely generated group $G \subset \psl$ with generating system
  $x_1$,...,$x_g$, the relation $[x_1,x_1]\cdots [x_g,x_g]=1$
  defines a representation $\rho:\Gamma_g\rightarrow\psl$ whose
  image is $G$, in the following way:
  $$\rho(a_i)=x_i\text{, }\rho(b_i)=x_i\text{, for }1\leq i\leq g$$
  In other words, every finitely generated subgroup
  of $\psl$ is the image of some representation of Euler class $0$.
\end{ps}

\begin{ps}
  Given two representations $\rho\in R_g$, $\rho'\in R_{g'}$, we have
  $[\rho(a_1),\rho(b_1)]\cdots[\rho(a_g),\rho(b_g)]=1$ and
  $[\rho(a'_1),\rho(b'_1)]\cdots[\rho(a'_{g'}),\rho(b'_{g'})]=1$
  so $[\rho(a_1),\rho(b_1)]\cdots[\rho(a'_{g'}),\rho(b'_{g'})]=1$.
  This defines naturally a representation in $R_{g+g'}$.
  This construction, together with the preceding remark, shows that
  there are also many subgroups of $\psl$ arising as the image of
  representations of non-zero Euler class, even in a fixed $R_g$.
\end{ps}

Now we prove proposition 1.1:
\begin{prop}
  Let $k$ be a fixed non-zero integer. There are only finitely many Fuchsian
  groups $\Gamma$ such that there exists a representation
  of $\Gamma_g$, for some $g\geq 2$,
  with image contained in $\Gamma$ and with Euler class $k$.
\end{prop}
\begin{proof}
  In other words, we have to show that 
  $e(\Gamma)\leq k$ holds true only for finitely many cocompact 
  Fuchsian groups $\Gamma$. We will prove that the inequality
  $0<d\left(2g-2+\displaystyle{\sum_{i=1}^r}
  \left(1-\frac{1}{k_i}\right)\right)\leq k$
  implies that $g$, $r$ and $d$ are bounded in terms of $k$. 
  Since $d\geq k_i$ for all
  $i\leq r$, every integer from the signature of $\Gamma$ will take only 
 finitely many values.
  
  First, $k_i\geq 2$ for all $i$ so
  $\displaystyle{\sum_{i=1}^r}\left(1-\frac{1}{k_i}\right)\geq\frac{r}{2}$,
  and $d\geq 1$ so $k\geq 2g-2+\frac{r}{2}$, i.e. $4g+r\leq 2k+4$. Thus
  $r$ and $g$ are bounded. We claim that $d\leq 42k$.
  
  \begin{enumerate}
    \item
    If $g\geq 1$ then $k\geq d\frac{r}{2}$. If $r=0$ then we have $d=1$,
    otherwise $k\geq\frac{d}{2}$ so $d\leq 42k$.
    \item
    Suppose $g=0$. If $d\geq 4k$, then
    $-2+\displaystyle{\sum_{i=1}^r}\left(1-\frac{1}{k_i}\right)\leq\frac{1}{4}$
    so $-2-\frac{r}{2}\leq\frac{1}{4}$, hence $r\leq 4$. But
    $-2+\displaystyle{\sum_{i=1}^r}\left(1-\frac{1}{k_i}\right)>0$ (see
    remark 4.3) so $r\geq 3$. Therefore, it suffices to consider 
    the cases $r=3$ and $r=4$.
  
    \begin{enumerate}
      \item
      If $r=4$, then
      $2-\frac{1}{k_1}-\frac{1}{k_2}-\frac{1}{k_3}-\frac{1}{k_4}>0$
      and $k_i\geq 2$. Hence, one of the $k_i$'s is greater than $3$.
      Therefore,
      $2-\frac{1}{k_1}-\frac{1}{k_2}-
      \frac{1}{k_3}-\frac{1}{k_4}\geq\frac{1}{6}$,
      and hence $d\leq 6k$.
      \item
      If $r=3$, then $1-\frac{1}{k_1}-\frac{1}{k_2}-\frac{1}{k_3}>0$.
      First suppose that
      $k_1$, $k_2$, $k_3\geq 3$. Then one of the $k_i$'s is greater than $4$,
      and hence $1-\frac{1}{k_1}-\frac{1}{k_2}-\frac{1}{k_3}\geq\frac{1}{12}$,
      so $d\leq 12k$.
      The only remaining case is when one of the $k_i$'s, say $k_3$,
      equals $2$.
      We then have
      $0<d\left(\frac{1}{2}-\frac{1}{k_1}-\frac{1}{k_2}\right)\leq k$.
      Thus, $k_1$, $k_2\geq 3$. And if $k_1$, $k_2\geq 4$, the condition
      $\frac{1}{2}-\frac{1}{k_1}-\frac{1}{k_2}>0$ forces one of $k_1$, $k_2$
      to be at least $5$, so that
      $\frac{1}{2}-\frac{1}{k_1}-\frac{1}{k_2}\geq\frac{1}{20}$ so $d\leq 20 k$.
      Hence, we may suppose that $k_1$ or $k_2$ (say, $k_2$)
      equals $3$. Now,
      $0<d\left(\frac{1}{6}-\frac{1}{k_1}\right)\leq k$ and hence
      $d\leq 42 k$.
    \end{enumerate}
  \end{enumerate}
\end{proof}

\begin{cor}
  No Fuchsian group can be injectively mapped
  into infinitely many Fuchsian groups.
  Moreover, if $\Gamma$ can be a sub-Fuchsian group of $\Gamma'$
  then $e(\Gamma')$ divides $e(\Gamma)$.
\end{cor}

\begin{ps}
  However, we can produce chains of inclusions of Fuchsian groups,
  of arbitrary length. For instance, one may check easily that
  $(2^k;-)\hookrightarrow (2^{k-1};2)\hookrightarrow (2^{k-2},4)
  \hookrightarrow\cdots\hookrightarrow (1;2^k)$.
\end{ps}

\subsection{Paucity of discrete representations}{\ }

\begin{prop}
  Let $g\geq 2$ and $k$ be such that $|k|\leq 2g-3$, $k\neq 0$.
  The representations
  of discrete image form a nowhere dense closed subset of $e^{-1}(k)$
  in $R_g$.
\end{prop}

\begin{proof}
  First, we show that the set of discrete representations is closed in
  $e^{-1}(k)$.
  For this, we recall J{\o}rgensen's lemma (see e.g. \cite{Katok}).
  If $A$, $B\in\psl$, we set
  $J(A,B)=|\Tr^2(A)-4|+|\Tr([A,B])-2|$.
  \begin{lem}
    If $\langle A,B\rangle$ is a non-elementary discrete group,
    then $J(A,B)\geq 1$.
  \end{lem}
  
  Now, let $\rho_n$ be a sequence
  of discrete representations of Euler class $k$, which converges pointwise
  to a representation $\rho$.
  Then $e(\rho)=k\neq 0$.
  Therefore, the representation $\rho$ is not elementary. Thus, there
  exist $T_1$, $T_2\in\rho(\Gamma_g)$ such that $[T_1,T_2]\neq 1$.
  Suppose that $\rho$ is not discrete. Then there exists
  a sequence $S_n\in\rho(\Gamma_g)$, $S_n\neq 1$, converging to $1$ in
  $\psl$. Hence, there exists $N\geq 0$ such that $J(S_N,T_1)<1$
  and $J(S_N,T_2)<1$. Now, $\rho_n$ converges to $\rho$
  pointwise so there exists
  $n_0\geq 0$ and $s,t_1,t_2\in\rho_{n_0}(\Gamma_g)$ such that
  $J(s,t_1)<1$, $J(s,t_2)<1$ and $[t_1,t_2]\neq 1$ (since these three
  conditions are open conditions). But this contradicts the assumption
  that $\rho_{n_0}$ has discrete image.

  Now we show that the closed set of discrete representations of Euler
  class $k$ is nowhere dense. For this, we show that it is the union
  of a countable family of closed and nowhere dense sets. Hence, by
  Baire's theorem, its complement is a dense set.

  Every discrete representation splits in the following way:
  $\Gamma_g\displaystyle{\mathop{\rightarrow}^h}\Gamma\hookrightarrow\psl$,
  where there are only finitely many groups $\Gamma$ possible.
  Now, for each Fuchsian group $\Gamma$,
  the set of morphisms $h:\Gamma_g\rightarrow\Gamma$ is countable. And
  the set of injective representations with discrete image
  $\Gamma\rightarrow\psl$ has an algebraic structure, similar to that
  of $R_g$.
  When it is pulled back by a morphism
  $h:\Gamma_g\rightarrow\Gamma$, we obtain an algebraic subset of $R_g$, of
  codimension at least $1$.
  It follows that for a fixed $h_0:\Gamma_g\rightarrow\Gamma$, the set
  of representations $\rho$ which factorize through $h_0$
  is a closed and nowhere dense subset of $R_g$.
  Now, the union of these sets, over the countable family of
  possible $\Gamma$ and $h$, is the set of discrete representations.
  As claimed, it is a countable union of nowhere dense subsets of
  $e^{-1}(k)\subset R_g$.
\end{proof}

\begin{ps}
  In the set $e^{-1}(0)$, the subset of representations with discrete
  image is not closed. Indeed, let $h:\Gamma_g\rightarrow\dZ^2$
  be a surjective mapping and let $q\in\dR$ be an  irrational number.
  The number $q$ can be approximated by rational numbers
  $\frac{\varphi_n}{\psi_n}$, with ${\rm gcd}(\varphi_n,\psi_n)=1$.
  Let $R(t)$ be a hyperbolic one-parameter subgroup of $\psl$. Define
  $\phi_n:\dZ^2\rightarrow\psl$ by $\phi_n(1,0)=R(1)$ and
  $\phi_n(0,1)=R(\frac{\varphi_n}{\psi_n})$. The representations
  $\phi_n\circ h:\Gamma_g\rightarrow\psl$ are elementary hence their
  Euler class is $0$. Moreover, they are all discrete. But they
  converge to a non-discrete representation.
\end{ps}

\subsection{Existence of discrete representations with prescribed 
Euler class}{\ }

Now, we prove theorem 1.3:
\begin{prop}
  For all $g\geq 2$ and $2-2g\leq k\leq 2g-2$, there exists a
  representation $\rho\in R_g$ with discrete image and Euler class $k$.
\end{prop}

\begin{proof}
  First, if $\rho$ is a representation of $\Gamma_g$ of Euler class $k$, then
  $\rho'(a_i)=\rho(b_{g-i})$, $\rho'(b_i)=\rho(a_{g-i})$ defines a
  representation of $\Gamma_g$ of Euler class $-k$, by Milnor's algorithm.
  Thus we just need to consider the case $k\geq 0$.
  
  If $k$ is even, set $k=2l$. Notice that $l+1\leq g$, by the Milnor-Wood
  inequality.
  Take a representation in the
  Teichm\"uller component of the surface of genus $l+1$.
  Its Euler class is $2(l+1)-2=k$. Now extend it by identity elements
  (as in remark 4.7)
  to get a representation of $\Gamma_g$. It still has discrete
  image, and by Milnor's algorithm its Euler class is $k$.
  
  Similarly, for an arbitrary $g$, it suffices  to find a representation
  of Euler class $2g-3$ in order to cover the case of odd $k$. We
  will extend it to representations of higher genus and the same 
Euler class, by using identity elements.
  \begin{enumerate}
  \item If $g$ is even, say, $g=2g'$, consider the Fuchsian group
 $\Gamma$  of signature $(g';2)$. We then have the following relations 
that hold in the lift $\widetilde{\Gamma}$ of $\Gamma$ to $\widetilde{\psl}$:
  $q^2=z$, $qc=z^{2g'-1}$, where $c$ is the product of  the $g'$ commutators
  defining the Fuchsian group $\Gamma$. 
Now, $z$ is central in $\widetilde{\psl}$
  so $cq=z^{2g'-1}$, and hence $cq^2 c=z^{4g'-2}$, i.e.
  $c^2=z^{4g'-3}=z^{2g-3}$. This relation enables us to write
  the element $z^{2g-3}$ using $2g'=g$ commutators in the
  Fuchsian group $\Gamma$ of signature $(g';2)$, 
and this defines a representation
  of $\Gamma_g$ of Euler class $2g-3$.
  Its image is a subgroup of $\Gamma$, hence discrete.
  Actually, its image is precisely $\Gamma$. 
\item   
  If $g$ is odd, say $g=2g'+1$, consider the Fuchsian group $\Gamma$ of
  signature $(g';2,2,2)$. In its lift $\widetilde{\Gamma}$ 
 in $\widetilde{\psl}$
  we have the following relations:
  $q_1^2=z$, $q_2^2=z$, $q_3^2=z$, $q_1 q_2 q_3 c=z^{2g'+1}$, where
  $c$ is once again a product of $g'$ commutators. Now, 
  $z$ is central so $q_2 q_3 c q_1=z^{2g'+1}$. Thus 
  $q_2 q_3 c q_2 q_3 c=z^{4g'+1}$ and hence
  $$z^{4g'-1}=q_2^{-1} q_3^{-1} c q_2 q_3 c=q_2^{-1}q_3^{-1}q_2 q_3
  q_3^{-1}q_2^{-1}c q_2 q_3 c=\left[q_2^{-1},q_3^{-1}\right]\cdot
  \left((q_2 q_3)^{-1} c (q_2 q_3)\right)\cdot c$$
  This implies that $z^{4g'-1}$ is the product of $2g'+1$ commutators, 
because the conjugation by $q_1 q_2$ enters the product of
  commutators $c$. Thus we obtain  
  a representation of $\Gamma_g$ of 
  Euler class $4g'-1=2g-3$. 
  Moreover, we can check that the
  image of this representation is the group $\Gamma$.
\end{enumerate}
\end{proof}

\subsection{Non-faithful representations are dense in every non-Teichm\"uller
component}{\ }

We will now prove theorem 1.1, which, for
the sake of completeness, we restate here:

\begin{prop}
  For all $g\geq 2$ and all $k$ such that $|k|<2g-2$,
  non-injective representations form a dense subset of the
  connected component $e^{-1}(k)$.
\end{prop}

\begin{proof}
  Suppose the contrary. Then there exists an open set $V\subset R_g$
  consisting only of
  injective representations of Euler class $k$ ($|k|< 2g-2$).
  Let $\rho_0\in V$. The representation $\rho_0$ is not discrete,
  so there exists $x\in\Gamma_g$ such
  that $\rho_0(x)\in Ell$
  (see e.g. \cite{Katok}). And $\rho_0$ is faithful so
  $\Tr(\rho_0(x))\in[0,2)=2\cos\theta_0$, with $\theta_0$ irrational
  (otherwise, $\rho_0(x)$ would be of finite order, which is impossible
  since $\Gamma_g$ is torsion-free).
  Now $\rho\mapsto\Tr(\rho(x))$ is continuous, so the angle $\theta$
  has to be constant on $V$ in order to keep being irrational.
  Hence the algebraic function $\Tr(\rho(x))$ is constant on the open
  set $V$ of the algebraic set $R_g$. Hence it is constant on the
  whole connected component $e^{-1}(k)$, since this connected component
  is contained within one irreducible component of the algebraic set $R_g$
  (see \cite{DBK}).
  
  It follows that for every $\rho\in e^{-1}(k)$, $\rho$ sends $x$
  to an elliptic element corresponding to an irrational angle. In
  particular, $\rho$ cannot have a discrete image. But this
  contradicts theorem 1.3.
\end{proof}

\begin{ps}
  As it was pointed out in \cite{DBK,Suoto}, injective representations form a
  dense subset of $R_g$.
  Therefore one may understand the set of non-faithful representations
  in $e^{-1}(k)$, for $|k|\leq 2g-1$ as $\dQ$ in $\dR$. A very
  comparable result has been recently proved by Glutsyuk in
  \cite{Glu}, for the set of representations of free groups into
  general Lie groups, thereby answering a question of Ghys.
\end{ps}

\section{Representations of odd Euler class}

The case of odd Euler class is somewhat simpler than the general case.
Even though the following results can be deduced from the latest part,
the approach is often more elementary.

\subsection{A cohomological criterion}{\ }

The commutator map
$$\begin{array}{ccc} \sl\times \sl & \rightarrow &
    \sl \\ A,B & \mapsto & [A,B] \end{array}  $$
defines an application
$$\begin{array}{ccc} \psl\times \psl & \rightarrow &
    \sl \\ A,B & \mapsto & [A,B] \end{array}  $$
Hence, there exists a continuous function
$\varepsilon : R_g\rightarrow\{-1,1\}$ such that for all
$\rho\in R_g$,
$$[\rho(a_1),\rho(b_1)]\cdots[\rho(a_g),\rho(b_g)]=
\varepsilon(\rho)$$

\begin{prop}
  For all $\rho\in\ho$, $\varepsilon(\rho)=(-1)^{e(\rho)}$, where
  $e(\rho)$ is the Euler class of $\rho$.
\end{prop}

\begin{proof}
  This follows directly from Milnor's algorithm: the element $z$
  maps to $-1$ in $\sl$.
\end{proof}

In particular, this gives a decomposition of $R_g$
into two sub-varieties of $(\sl)^{2g}$. One of them is defined by the
equation $[x_1,y_1]\cdots[x_g,y_g]=1$ and the other is defined by
$[x_1,y_1]\cdots[x_g,y_g]=-1$. These two algebraic varieties are
irreducible.
In fact the invariant $\varepsilon(\rho)$ is precisely the second
Stiefel-Whitney class $w_2(\rho)$ (see \cite{Gold84,Goldman,DBK}).

\begin{cor}
  Let $\Gamma$ be a subgroup of $\psl$ of finite type and let
  $\widetilde{\Gamma}$ be its lift to $\sl$ containing
  $-1$. (Note that this is well-defined).
  Then the two following assertions are equivalent :
  \begin{itemize}
    \item there exists $\rho\in R_g$ of odd Euler class and whose image
    is $\Gamma$
    \item $-1$ is a product of commutators in $\widetilde{\Gamma}$, i.e.
    $-1$ maps to the neutral element in $H^1(\widetilde{\Gamma})$.
  \end{itemize}
\end{cor}

\begin{proof}
  If there exists a representation $\rho\in R_g$ of odd Euler class
  whose image is $\Gamma$ then by proposition 4.1,
  $[\rho(a_1),\rho(b_1)]\cdots[\rho(a_g),\rho(b_g)]=-1$ so $-1$ is
  indeed a product of commutators in the image of $\rho$.

  Conversely, suppose $-1$ is a product of commutators in the
  group $\Gamma$: $[x_1,y_1]\cdots[x_g,y_g]=-1$. Let
  $z_1$,..., $z_n$ be a system of generators of $\Gamma$.
  Then
  $[x_1,y_1]\cdots[x_g,y_g][z_1,z_1]\cdots[z_n,z_n]$ defines
  a representation $\rho$ of $\Sigma_{g+n}$, whose image is exactly $\Gamma$.
\end{proof}

\begin{cor}
  There are no representations of odd Euler class whose image
  is in $PSL(2,\dZ)$.
\end{cor}

\begin{proof}
  The element $-1\in SL(2,\dZ)$ maps to a non-trivial element in
  the group $H^1(SL(2,\dZ))\cong\dZ/12\dZ$.
\end{proof}

\subsection{An explicit example}{\ }

We can use this characterization to give an explicit example of
a discrete representation of odd Euler class as follows.
  Let $\Gamma$ be the Fuchsian triangular group of signature 
$(0;2,3,7)$. 
  Then $\Gamma$ has the following presentation:
$\Gamma=\left\langle q_1,q_2,q_3\left|q_1^2=1,q_2^3=1,q_3^7=1,q_1q_2q_3=1
\right.\right\rangle$.
  Its lift $\widetilde{\Gamma}$ then has the following presentation:
$$\widetilde{\Gamma}=\left\langle q_1,q_2,q_3,h\left|
\begin{array}{l}h^2=1,h q_1=q_1 h,h q_2=q_2 h,
h q_3=q_3 h, \\ q_1^2=h^{\beta_1},
q_2^3=h^{\beta_2},q_3^7=h^{\beta_3},q_1q_2q_3=h^{\beta}
\end{array}\right.\right\rangle$$
  One can replace $q_2$ and $q_3$, respectively, by
$q_2 h$ and $q_3 h$, to get $\beta_2=\beta_3=1$,
and then replace $q_1$ by $q_1 h$ to get $\beta=1$ in this
presentation.
  Now, if we had $\beta_1=0$, this would mean that
$q_1$ has order $2$ in $\sl$ and thus is equal to $1$ or $-1$,
and hence is $1$ in $\psl$, which is not the case. Thus $\beta_1=1$
(mod 2). Finally, we get the following presentation of
$\widetilde{\Gamma}$:
$$\widetilde{\Gamma}=
\left\langle q_1,q_2,q_3,h\left|
\begin{array}{l}h^2=1,h q_1=q_1 h,h q_2=q_2 h,
h q_3=q_3 h, \\ q_1^2=h,q_2^3=1,q_3^7=1,q_1q_2q_3=1\end{array}
\right.\right\rangle$$
  Consider now $H^1(\widetilde{\Gamma})$. In this group,
$q_1 q_2 q_3=1$ so
$(q_1 q_2 q_3)^{42}=1=
q_1^{42}q_2^{42}q_3^{42}=h^{21} q_2^{3\times 14} q_3^{7\times 6}=
h.h^{2\times 10}=h$
so $h=1$ in $H^1(\widetilde{\Gamma})$. Equivalently,
$-1$ is a product of commutators in $\widetilde{\Gamma}$.
In the same way as in remark 4.3, we can explicitly write
$-1$ as a product of commutators, thus defining a representation of
some $\Gamma_g$ in $\psl$.

\subsection{Discrete representations of odd Euler class}{\ }

If $\Gamma$ is a cocompact Fuchsian group with signature
$(g;k_1,k_2,\ldots,k_r)$, let $m(\Gamma)$ be the maximal power
of $2$ dividing one of the $k_i$'s. If $m(\Gamma)=0$, set
$n(\Gamma)=0$. Otherwise, let $n(\Gamma)$
be the number
of $k_i$'s which are divisible by $2^{m(\Gamma)}$.
We then have the following characterization:

\begin{prop}
  A Fuchsian group $\Gamma \subset \psl$ is the image of a
  representation of $\Gamma_g$, for some $g$, if and only if $\Gamma$ is
  a cocompact Fuchsian group such that $n(\Gamma)$ is odd.
\end{prop}

\begin{proof}
  This may be deduced directly from proposition 4.1, but we give a
  slightly different proof here. First, by remark 4.4
  (or proposition 4.1), $\Gamma$ has
  to be cocompact (otherwise $H^2(\Gamma)=0$ and $e(\rho)=0$).
  
  Now, let $\Gamma$ be a cocompact Fuchsian group, with signature
  $(g;k_1,k_2,\ldots,k_r)$. This means that $\Gamma$ has the following
  presentation:
  $$\Gamma=\left\langle q_1,q_2,\ldots,q_r,a_1,\ldots,b_g
  \left|q_1^{k_1},
  \ldots,q_r^{k_r},q_1\cdots
  q_r[a_1,b_1]\cdots[a_g,b_g]\right.\right\rangle$$

  Now using lemma 4.1 and the fact that $h^2=1$ here (we are in
  $\sl$ and not in $\widetilde{\psl}$), we get the following presentation
  for the lift $\widetilde{\Gamma}$ in $\sl$:
  $$\widetilde{\Gamma}=\left\langle q_1,\ldots,q_r,a_1,
  \ldots,b_g,h\left|
  \begin{array}{l}hq_1h^{-1}q_1^{-1},\ldots,
  hb_gh^{-1}b_g^{-1},h^2,q_1^{k_1}h,\ldots,q_r^{k_r}h,\\
  q_1\cdots q_r[a_1,b_1]\cdots[a_g,b_g]h^{r}\end{array}\right.\right\rangle$$
  It follows that the abelianisation of $\widetilde{\Gamma}$ has the following
  abelian presentation:
  $$H^1(\widetilde{\Gamma})=
  \left\langle q_1,\ldots,q_r,h \left| h^2,q_1^{k_1}h,\ldots,q_r^{k_r}h,
  q_1\cdots q_r h^r\right.\right\rangle^{\rm ab}$$
  Now reorder the $q_i$'s so that the powers of $2$ dividing $k_i$ are
  decreasing.
  For $1\leq i\leq r$, let $k_i=2^{u_i}v_i$ with $v_i$
  odd. In particular, $ u_i=m(\Gamma)$ if $1\leq i\leq n(\Gamma)$.

\begin{enumerate}
 \item  First suppose that $n(\Gamma)$ is odd. Then
  $\left(q_1\cdots q_r h^r\right)^{2^{m(\Gamma)}v_1
  \cdots v_r}=h^{n(\Gamma)}=1$ so $h=1$ in
  $H^1(\widetilde{\Gamma})$.
\item 
  Assume now  that $n(\Gamma)$ is even. Then we will define
  a morphism
  $\phi:H^1(\widetilde{\Gamma})\rightarrow \dS^1=\{z\in\dC:|z|=1\}$.
  First, let $\phi(h)=-1$. Next, if $i\geq 2$, let
  $\phi(q_i)=\expon{\frac{i\pi}{2^{u_i}}}$. Then
  $\phi(q_i^{k_i}h)=-\expon{\frac{i\pi k_i}{2^{u_i}}}=\expon{i\pi(1+v_i)}=1$.
  If $n(\Gamma)=0$ (or equivalently, if all the $k_i$'s are odd),
  set $\phi(q_1)=-1$. Then one may easily check that $\phi(q_1^{k_1}h)=1$
  and $\phi(q_1\cdots q_r h^r)=1$,
  so our morphism $\phi$ is well-defined.
  Otherwise, let
  $\phi(q_1)=\expon{-\frac{i\pi}{2^m}((n-1)+2^m r+\sum_{i=n+1}^{r} 2^{m-u_i})}$.
  Then once again we can check that $\phi(q_1^{k_1}h)=1$ and that
  $\phi(q_1\cdots q_r h^r)=1$, so our morphism is well-defined.
  This proves that $h\neq 1$ in $H^1(\widetilde{\Gamma})$.
\end{enumerate}
\end{proof}

\begin{ps}
  As in remark 4.3, in the case $h=1$ in $H^1(\widetilde{\Gamma})$,
  we can define an explicit representation of odd Euler class, whose
  image is $\Gamma$. This is done by writing
  $(q_1\cdots q_r)^d q_1^{-d}\cdots q_r^{-d}$,
  with $d=2^{m(\Gamma)}v_1\cdots v_r$, as a product of
  commutators.
\end{ps}


\begin{thebibliography}{Name YY}
  \bibitem{CullerShallen}M. Culler, P.B. Shallen: Varieties
of group representations and splittings of $3$-manifolds,
{\it Ann. of Math.},
{\bf 117}(1983), 109-146.
  \bibitem{DBK}J. DeBlois, R. Kent IV: Surface groups are
frequently faithful, {math.GT/0411270}, 2004.
  \bibitem{Brown}K. Brown: {\it Cohomology of groups}, GTM 87, 
Springer-Verlag, 1994.
  \bibitem{Galloetal}D. Gallo, M. Kapovich, A. Marden: The monodromy
groups of Schwarzian equations on closed Riemann surfaces, {\it Ann. of
Math.} {\bf 151}(2000),  625-704.
  \bibitem{Ghys}E. Ghys: Classe d'Euler et minimale exceptionnel,
{\it Topology} {\bf 26}(1987), 93-105.
  \bibitem{Glu}A. Glutsyuk: Instability of nondiscrete free subgroups in Lie groups, {math.DS/0409556}, 2004.
  \bibitem{Gold3}W. Goldman: {\it Discontinuous groups and the Euler
class}, doctoral thesis, Berkeley, 1980.
  \bibitem{Gold84}W. Goldman: The symplectic nature of fundamental
groups of surfaces, {\it Advances  Math} {\bf 54}(1984), 200-225.
  \bibitem{Goldman}W. Goldman: Topological components of spaces
of representations, {\it Invent. Math.} {\bf 93}(1988), 557-607.
  \bibitem{Hitchin}N.J. Hitchin: The self-duality equations on a
Riemann surface, {\it Proc. London Math. Soc.} (3) {\bf 55}(1987), 59-126.
  \bibitem{Katok}S. Katok: {\it Fuchsian groups}, Chicago Lectures Math., 
1992.
  \bibitem{Magnus}W. Magnus: Rational representations of Fuchsian
groups and non-parabolic subgroups of the modular group, {\it Nachr. Akad. Wiss.
Gottingen Math.-Phys. K.I. II}, 1973, 179-189.
  \bibitem{Milnor}J.W. Milnor: On the existence of a connection
with curvature zero, {\it Comment. Math. Helv.} {\bf 32}(1958), 215-223.
  \bibitem{Suoto}E. Breuillard, T. Gelander, J. Souto, P. Storm: Dense subgroups of Lie groups, {\it in preparation}.
  \bibitem{Tan}S.P. Tan: Branched $CP^1$-structures on surfaces with
prescribed real holonomy, {\it Math. Ann.} {\bf 300}(1994), 649-667.
  \bibitem{Wood}J.W. Wood: Bundles with totally disconnected
structure group, {\it Comment. Math. Helv.} {\bf 51}(1971), 183-199.
\end{thebibliography}
\end{document}